\documentclass[graybox]{svmult}

\usepackage[export]{adjustbox}

\usepackage[round]{natbib}
\usepackage[utf8x]{inputenc}
\usepackage{amsmath, amssymb,bbm,wasysym}
\usepackage{arydshln}

\usepackage{mathptmx}       
\usepackage{helvet}         
\usepackage{courier}        
\usepackage{type1cm}        
\usepackage{makeidx}         
\usepackage{graphicx}        
\usepackage{multicol}        
\usepackage[bottom]{footmisc}

\newcommand{\ZZ}{{\mathbbm{Z}}} \newcommand{\NN}{{\mathbbm{N}}}

\makeindex             

\begin{document}

\title*{Magic squares of subtraction of Adam Adamandy Kocha\'nski}
\author{Henryk Fuk\'s}
\institute{Henryk Fuk\'s \at Department of Mathematics and Statistics, Brock University,  St. Catharines,
ON, Canada \email{hfuks@brocku.ca}}
%
%
\maketitle

\abstract{The problem of the construction of magic squares occupied many mathematicians of the 17th century.
The Polish Jesuit and polymath Adam Adamandy Kocha{\'n}\-ski studied this subject too, and in 1686 he published 
a paper in \emph{Acta Eruditorum} titled ``Considerationes quaedam circa Quadrata et Cubos Magicos''. In that paper
he proposed a novel type of magic square, where in every row, column and diagonal, if the entries are sorted in decreasing order,
 the difference between the sum  of  entries with odd indices and those with even indices is constant.
He called them \emph{quadrata subtractionis}, meaning squares of subtraction. He gave examples of such squares of orders 4
and 5, and challenged readers to produce an example of square of order 6. We discuss the likely method which he used to produce
squares of order 5, and show that it can be generalized to arbitrary odd orders. We also show how to construct doubly-even squares.
At the end, we show an example of a square of order 6, sought by Kocha\'nski, and discuss the enumeration of squares of subtraction.
 }
\section{Introduction}

Contributions of Jesuits to science and mathematics have attracted considerable attention in recent years \citep{feingold2013,o2016jesuits}. It is now well documented that the Society of Jesus played an especially important role in the  intellectual life of Polish–Lithuanian Commonwealth  \citep{stasiewicz2004}, producing a number of notable writers, preachers,
theologians, missionaries, and scientists.  The Polish Jesuit Adam Adamandy Kocha\'nski (1631--1700) was one of 
such remarkable individuals. He pursued problems in a multitude  of diverse fields, including 
philosophy, mathematics, astronomy, philology, design and construction of mechanical clocks and computing devices, and many others.  He published relatively little, and  his published works remain rather unknown today,
even among historians of mathematics and science. In recent years, however, his accomplishments seem to have attracted
some more interest, in large part due to the efforts of B. Lisiak SJ, who reprinted all the published works of
 \citet{KochanskiOpuscula1,KochanskiOpuscula2,KochanskiOpuscula3} and who authored a comprehensive account of Kocha\'nski's life and work \citep{Lisiak05}. Moreover, B. Lisiak SJ and L. Grzebie\'n SJ  
gathered and published all the surviving letters of Kocha\'nski, including
correspondence with Gottfried Leibniz, Athanasius Kircher SJ, Johannes  Hevelius, Gottfried  Kirch, and many other luminaries of the 17th century \citep{LisiakGrzebien05}.

The most important and original mathematical works of Kocha\'nski appeared in \emph{Acta Eruditorum} between 1682 and 1696. Among those, three
are particularly interesting, namely papers published in \emph{Acta} in  \citeyear{KochanskiSolutio}, \citeyear{Kochanski1685}, and \citeyear{KochanskiConsiderationes} \citetext{all three reprinted in  \citealp{KochanskiOpuscula1}, English translation in \citealp{paper56}}. The 1685 paper,  \emph{Observationes cyclometricae}, is relatively well-known and is
often considered to be the most interesting one. It proposes for the rectification of the circle an approximate procedure 
which received some attention from both contemporaries of Kocha\'nski and  historians of mathematics \citep{Montucla1755,Cantor1880,Gunther1921}. What is generally less known, however, is that the same paper included an interesting sequence of rational approximations
of $\pi$, somewhat similar to continued fractions. Detailed discussion of this sequence and its properties can be found in a recent article
\citep{paper42}, so we will not consider it  here. We will only mention that the sequence of integers on which the aforementioned approximation of $\pi$ is based
is now  included in the Online Encyclopedia of Integer  Sequences as A191642 \citep{A191642}. I proposed to call it \emph{Kocha\'nski's sequence}.

The subject of this note is the 1686 paper, \emph{Considerationes quaedam circa Quadrata \& Cubos Magicos}. This work, like all of Kocha\'nski's
other papers, touches upon more than one subject.
In the first part, ``classical'' magic squares are discussed, and Kocha\'nski presents some previously unknown magic squares using the method of their construction  developed by A. Kircher. The second part is more original
and interesting,  because Kocha\'nski introduces there a new type of magic square, to be called \emph{quadrata subtractionis} or \emph{squares of subtraction}. An $n\times n$ 
square of subtraction is an arrangement of consecutive integers from 1 to $n^2$ in such a way that in every row, column and diagonal, if the entries are sorted in decreasing order,
 the difference between the sum  of  entries with odd indices and the sum of entries with even indices is constant. 
Kocha\'nski gives examples of $n=4$ and $n=5$ squares of subtraction, and then challenges fellow mathematicians to
produce an example of a square of order 6. To my knowledge, nobody has ever taken  this challenge, and nobody has ever studied magic squares of subtraction in the 330 years following the publication of \emph{Considerationes}.
The only author mentioning magic squares of subtraction in a published work is Z. Pawlikowska (\citeyear{Pawlikowska69}), who briefly discussed them in a paper describing the mathematical works of Kocha\'nski and who guessed the most likely method used by him to produce squares of order 5. 

The purpose of this article, therefore, is to bring magic squares of subtraction back to life, with the hope 
of stimulating others to study their properties and methods of construction. We will take a closer look at the
examples supplied by Kocha\'nski and present the general method of construction of squares of odd order and doubly-even order,
based on these examples. We will also give an example of an $n=6$ square, which he requested, and discuss the enumeration of squares of subtraction. 

\section{Notation and examples}
\begin{figure}
\begin{center}
 \includegraphics[width=7cm]{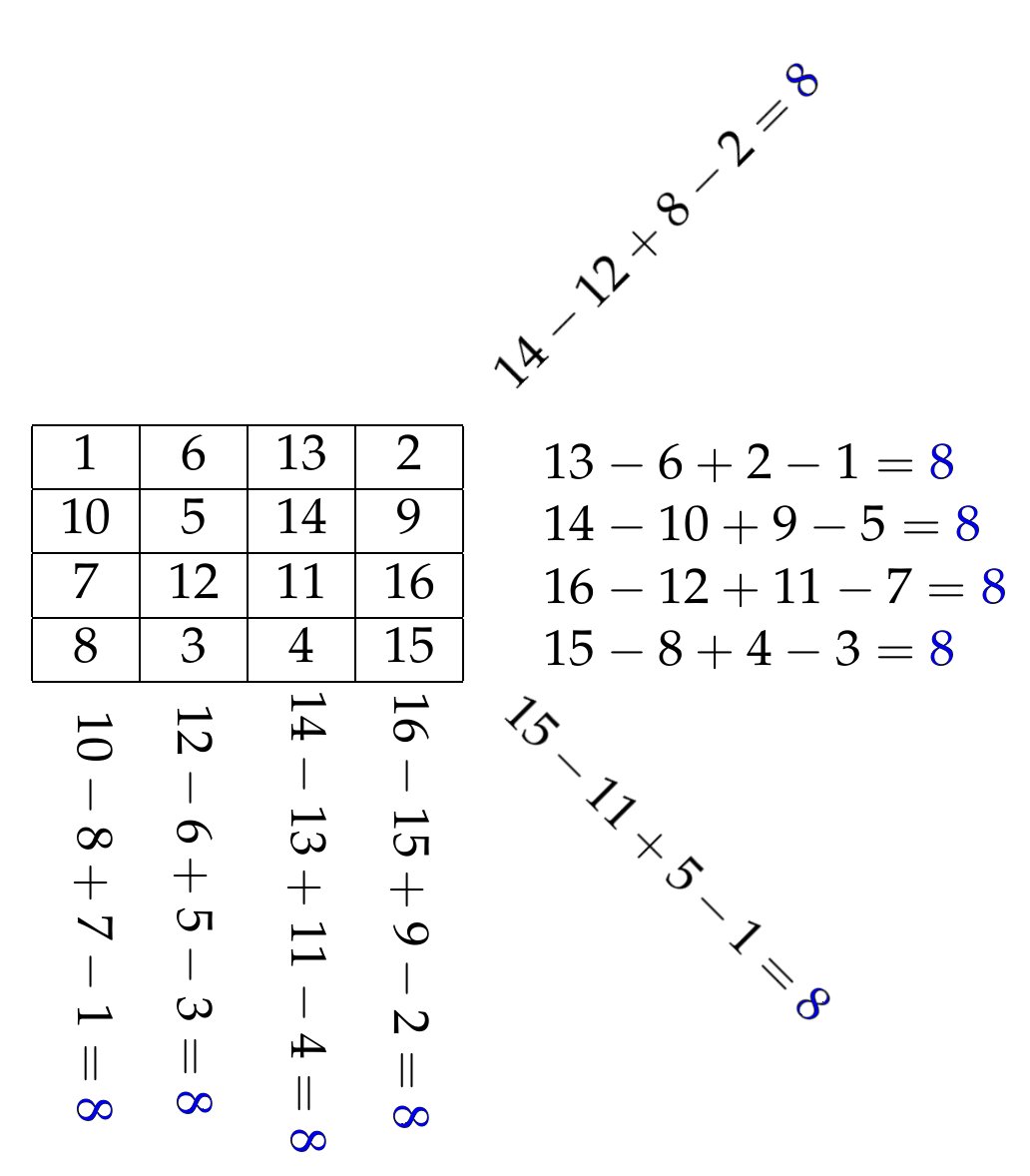}
\end{center}
\caption{Example of a square of subtraction of order $4$.}\label{exorder4}
\end{figure}

\begin{figure}
\begin{center}
 \includegraphics[width=11cm]{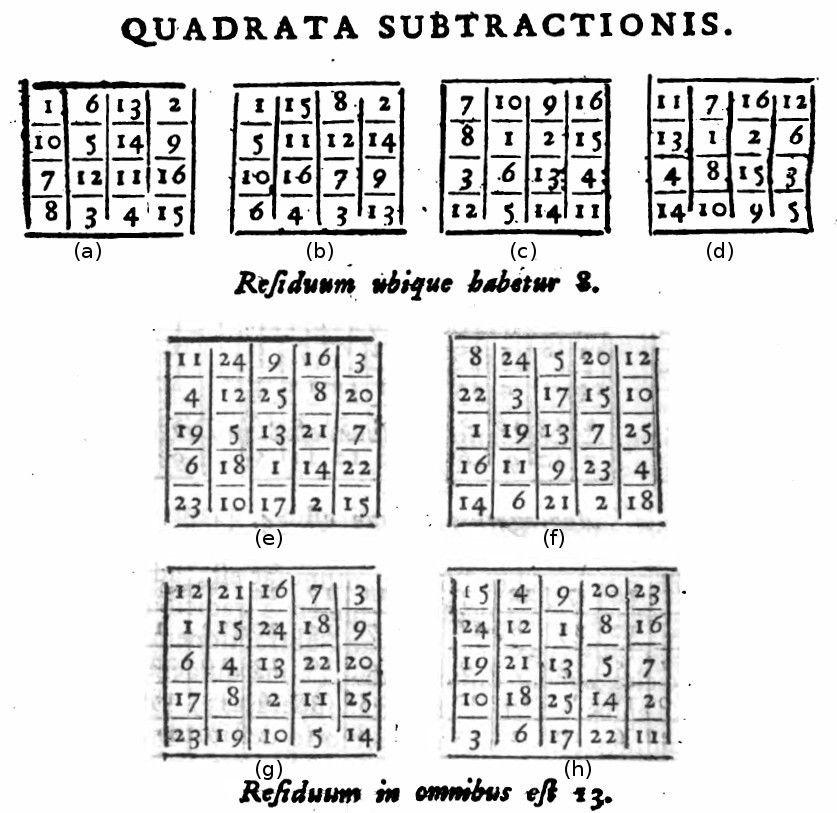}
\end{center}
\caption{Original examples of squares of orders 4 and 5 from Kocha\'nski's \emph{Considerationes}.} \label{originalexamples}
\end{figure}

Kocha\'nski defines  $n \times n$ squares of subtraction as arrangements of the consecutive integers $1,2, \ldots, n^2$ which have the same \emph{residuum}
in rows, columns, and diagonals. His definition of ``residuum'' is as follows.
\begin{quotation}
\noindent\emph{In his Quadratis Subtractione tractandis ita proceditur. Numerus in assumpta Columna, Trabe, vel Diagono minimus, subducendus est a proxime majore, residuum majoris hujus detrahatur a proxime consequente, atque ita porro : ultimum enim residuum est illud universale [...].}\\
In these Squares produced by Subtraction one proceeds as follows. The smallest number in a selected Column, Row, or Diagonal
is subtracted from the next-largest number, the difference is subtracted form the subsequent larger number, and so on: the 
final result is the universal number [...].
 \end{quotation}
This means, for example,  that if numbers in a given row, column or diagonal of a square $5 \times 5$,  sorted in decreasing order, are $a_1,a_2,a_3,a_4, a_5$, then the aforementioned residuum is
\begin{equation*}
 a_1-(a_2-(a_3-(a_4-a_5))).
\end{equation*}
Obviously, this can be written in a more convenient way, without nested parentheses.
Let us, therefore,  define the residuum of a vector $\mathbf{x}=(x_1, x_2, \ldots x_n)$ using modern notation.
First, we sort the vector $\mathbf{x}$ in decreasing order, obtaining  $\mathbf{\tilde{x}}$ such that 
 $\tilde{x}_1 \geq \tilde{x}_2 \geq \ldots \geq \tilde{x}_n$; then the
 \emph{residuum} of $\mathbf{x}$ is defined as
\begin{equation}
 \mathrm{res} (\mathbf{x})= \tilde{x}_1-\tilde{x}_2+\tilde{x}_3 - \ldots \tilde{x}_n= \sum_{i=1}^{n} (-1)^{i+1} \tilde{x}_i.
\end{equation} 
\emph{A Magic square of subtraction} is an arrangement of numbers  in a square array such that all rows, columns, and both diagonals have the same residuum. 
If entries of the array belong to the set $\{1,2,\ldots, n^2\}$, then the magic square of subtraction
will be called \emph{normal}. In what follows we will consider only normal squares (unless indicated otherwise), so we will omit the designation ``normal'' and we will refer to them simply as \emph{squares of subtraction}. 
The number $n$ will be called the \emph{order} of the square of 
subtraction. Figure \ref{exorder4} shows an example of a square of subtraction of order 4 with
residuum equal to 8.
This is actually one of the squares given by Kocha\'nski. In the original paper he gave examples of seven other squares of subtraction, of orders 4 and 5,
as shown in Figure~\ref{originalexamples}.
It is worth noticing that the smallest order he considered is 4, and that he does not mention any smaller orders.
He was likely aware that squares of order 2 and 3 do not exist, and this is not too hard to prove.

%
\section{Formal construction for odd $n$}
\begin{figure}
\begin{center}
 \includegraphics[width=6cm]{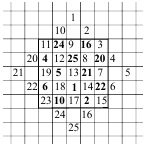}
\end{center}
\caption{Construction of a square of subtraction of order $5$. Relocated entries are shown in boldface.} \label{exconst}
\end{figure}

We will start with odd-order squares, as these are easier to construct.
Kocha\'nski gives examples of squares of order 5, and this is indeed  the smallest possible odd order. The first of these,   shown in Figure~\ref{originalexamples}e, is reproduced again in Figure~\ref{exconst},
as the framed square. As observed by \citet{Pawlikowska69}, it can be constructed by writing consecutive
numbers from 1 to 25 along diagonal lines of length 5, one below the other, with alternating direction. Afterwards, one needs to relocate the numbers which are outside the framed square to its interior,
by ``wrapping'' them back to the square (boldface numbers in Figure~\ref{exconst} are those which have been relocated).
\begin{figure}
 \begin{center}
 \includegraphics[width=6cm]{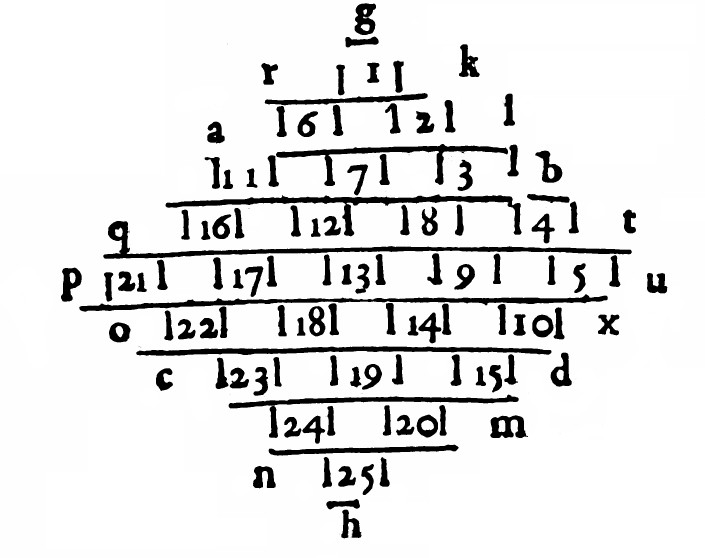}  \includegraphics[width=5cm]{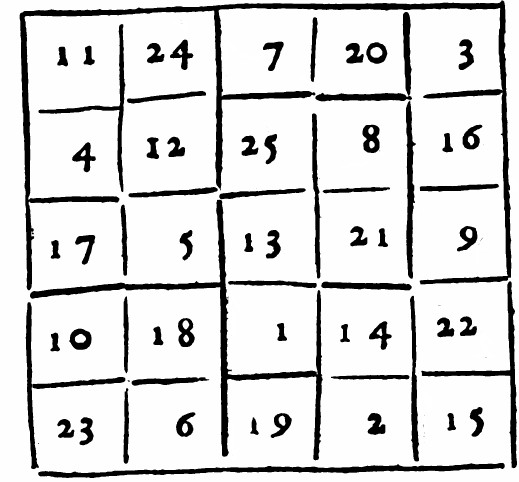}
\caption{Construction of a summation magic square of order 5 reproduced from pp. 89 and 90
  of Kircher's \emph{Arithmologia} \citep{kircher1665}. } \label{kircher}
\end{center}
\end{figure}

How did Kocha\'nski discover this method of construction? From his surviving letters we know that he was familiar with works of Athanasius Kircher SJ,
whom he greatly esteemed and with whom he maintained correspondence \citep{LisiakGrzebien05, Lisiak05}. In Kircher's book \emph{Arithmologia} \citep{kircher1665},
a method of construction of magic squares of summation is discussed, here reproduced in Figure~\ref{kircher}. In this method,
we write consecutive integers along diagonal lines, one below the other, but maintaining the same direction. Numbers
which fall outside of the square $n\times n$ are relocated to the interior of the square. This method has
been originally proposed by Claude Gaspard Bachet de Méziriac (\citeyear{Bachet1624}), but since Kocha\'nski does
not mention Bachet's name in \emph{Considerationes}, it seems likely that  he learned about it from Kircher (Kircher's
name is explicitly listed in the first paragraph of \emph{Considerationes} as one of  the ``ingenious men'' who have
studied magic squares).

Kocha\'nski noticed that by reversing the order in which we write numbers on each diagonal in Kircher's example, instead of a magic square of summation we obtain a square of subtraction. He did not prove this fact, but merely produced examples. 
We will now show that this method indeed produces the desired magic squares of subtraction, and that it can be generalized
to any odd order $n$. To simplify notation, we will define 
\begin{equation}
   [ k ]_n = \Big( (k- 1) \mod n \Big) +1. 
\end{equation}
Note that if $k \in \ZZ$ then $[k]_n\in \{1,2,\ldots,n\}$. This will be used to bring indices of matrices back to the range $\{1,2,\ldots,n\}$.

\begin{proposition}
 Let $n$ be an odd integer, $n\geq 3$, $m=(n+1)/2$, and let 
\begin{equation} \label{deff}
 f(i,j)=(m-i)(-1)^j + n (n-j) +m.
\end{equation}
Then the matrix 
\begin{equation}\label{defm}
M_{i,j}= \begin{cases}
           f\left( \frac{1}{2}(i+j),\frac{1}{2} (j-i) +m \right) & \text{if $i,j$ are of the same parity,} \\
           f\left( \left[\frac{1}{2}(i+j+n)\right]_n,
                      \left[\frac{1}{2} (j-i+n) +m\right]_n  \right) & \text{otherwise}.
          \end{cases}
\end{equation}
is a magic square of subtraction of order $n$ with residuum $\frac{n^2+1}{2}$.
\end{proposition}
We will be using Figure \ref{exconst} to illustrate the proof. 
Let us first 
explain where the formulae in eqs. (\ref{deff}) and (\ref{defm}) come from. 
Note that the magic square
of order 5 shown in Figure~\ref{exconst} has been obtained by rotating the matrix
\begin{equation} \label{defA}
A=\left(
\begin{array}{rrrrr} 
21 & 20 & 11 & 10 & 1\\ 
22 & 19 & 12 & 9 & 2\\ 
23 & 18 & 13 & 8 & 3\\ 
24 & 17 & 14 & 7 & 4\\  
25 & 16 & 15 & 6 & 5 \\ 
\end{array} \right)
\end{equation}
by 45 degrees counterclockwise. The function $f(i,j)$ defined in the proposition is
simply the indexing function used to construct $A$, so that $A_{i,j}=f(i,j)$. 
After $A$ is rotated, one needs to ``squeeze'' all entries of this matrix into
a square box $5 \times 5$, shown in Figure~\ref{exconst}. ``Squeezing'' is achieved
by relocating all entries which are outside the box to the inside of the box, by 
bringing their indices to the range $\{1,2,\ldots,5\}$ with $[\cdot]$ operator. 	
This simply means that each entry located outside the box at $(i,j)$ is relocated to
$([i]_n,[j]_n)$.  The result is the matrix $M$,
\begin{equation}
M=\left(
\begin{array}{rrrrr} 
11 & 24 & 9 & 16 & 3\\ 
4 & 12 & 25 & 8 & 20\\ 
19 & 5 & 13 & 21 & 7\\ 
6 & 18 & 1 & 14 & 22\\  
23 & 10 & 17 & 2 & 15 \\ 
\end{array} \right),
 \end{equation}
which is the desired magic square of subtraction. The rotation by 45 degrees 
around the origin can be described algebraically as the transformation $(i,j) \to
\left(\frac{1}{2}(i+j), \frac{1}{2}(j-i)\right)$, and we need to add $m$ to the second index to bring the origin of coordinate system to the right location, thus
\begin{equation}
M_{i,j} =f \left( \frac{1}{2}(i+j),\frac{1}{2} (j-i) +m \right) .
\end{equation}
Note that the above works only for entries which ended up inside the box, not needing 
relocation. These happen to be entries for which $i,j$ are both odd or both even.
The remaining entries must be relocated, and this is achieved by combined translation and
use of the operator $[\cdot]_n$, resulting in the second line of eq. (\ref{defm}).
\begin{figure}
\begin{center}
 \includegraphics[width=6cm]{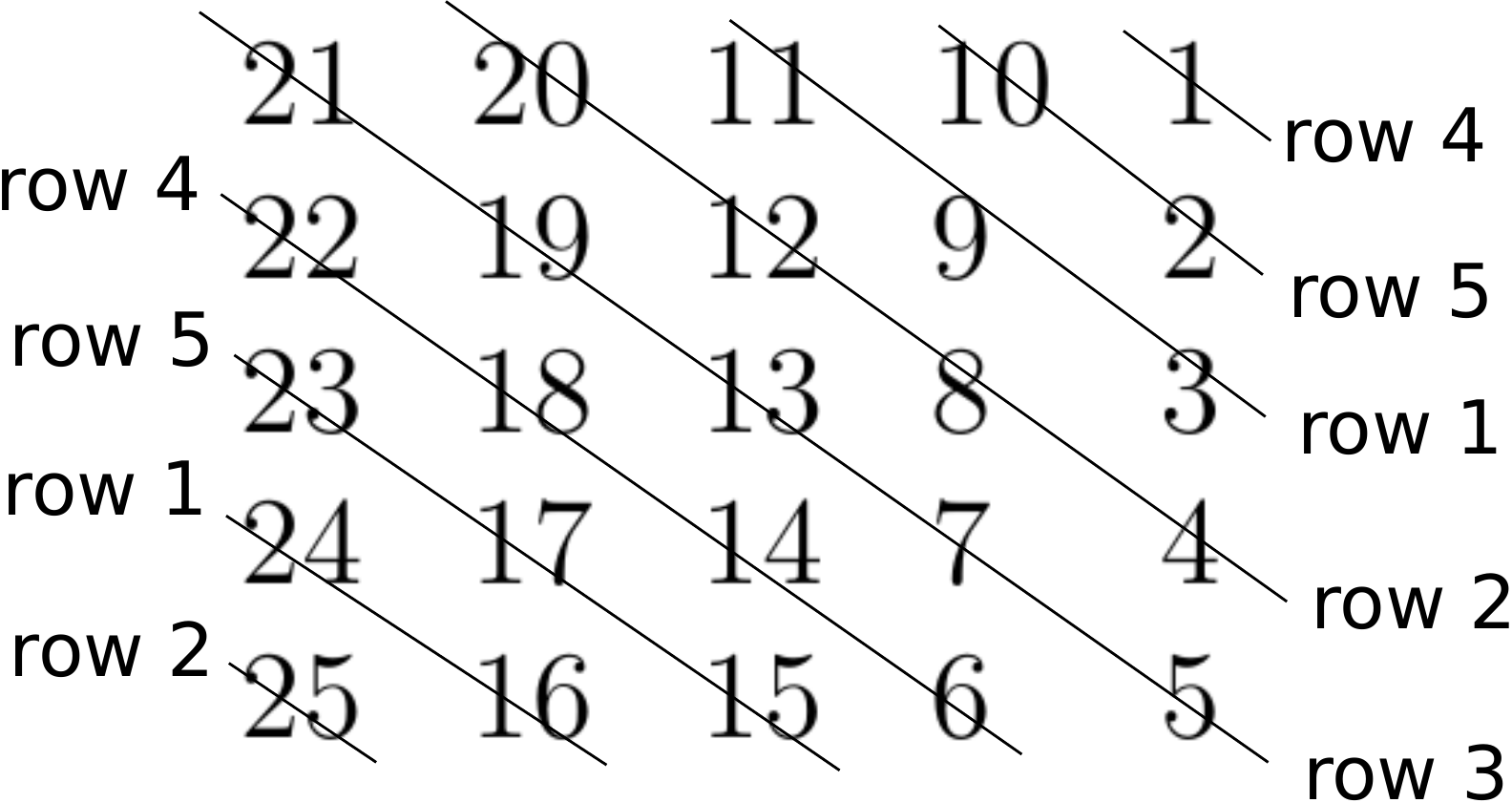}
\end{center}
\caption{Correspondence between diagonals of $A$ and rows of $M$.}\label{rowsassign}
\end{figure}

From Figure \ref{exconst}, one can see that the third row of $M$ is the diagonal of
$A$. Other rows of $M$ can be constructed from numbers lying on two diagonal segments
of $A$, as  illustrated in Figure \ref{rowsassign}.

This can be easily generalized to arbitrary $n$. For convenience, we will introduce 
index $p \in \{0,1, \ldots n-1\}$ labeling rows of $M$, such that the actual row number 
is equal $[p-m-1]_n$. The middle row of $M$ will then always correspond to $p=0$.
The set 
\begin{equation}
 R_p=\{ f(p+i,i)\}_{i=1}^{n-p} \cup \{ f(p+i-n,i)\}_{i=n-p+1}^{n},
\end{equation}
composed of two diagonal lines of $A$,  contains all entries of the row $[p-m-1]_n$ of $M$. Note
that the second component of the above set is empty when $p=0$, as in this case we only
need one diagonal line (the main diagonal).

It is important to notice that the elements of the set $R_p$ are already sorted in descending
order, thus we can easily compute its residuum,
\begin{equation}
 \mathrm{res} (R_p)=\sum_{i=1}^{n-p} (-1)^{i+1} f(p+i,i) 
+\sum_{i=n-p+1}^{n}  (-1)^{i+1} f(p+i-n,i)
\end{equation}
Using the definition of $f$ (eq. \ref{deff}) and the fact that $n$ is odd, one can easily compute
the above sums, obtaining
\begin{equation}
 \mathrm{res} (R_p)= \frac{n^2+1}{2}.
\end{equation}
This means that the residua of all rows of $M$ are the same and equal to $\frac{n^2+1}{2}$.

Exactly the same observations can be made for columns of $M$. We label them with the
integer $q \in \{0,1,\ldots, n-1\}$, so that the actual column number is $[3-q]_n$.
Then all the entries of the column $[3-q]_n$ of $M$ are contained in the set
\begin{equation}
 C_q=\{ f(n-i-q+1,i)\}_{i=1}^{n-q} \cup \{ f(2n-i+q+1,i)\}_{i=n-q+1}^{n}.
\end{equation}
As before, the above set is already sorted in descending order,
thus 
\begin{equation}
 \mathrm{res} (C_q)=\sum_{i=1}^{n-q} (-1)^{i+1} f(n-i-q+1,i) 
+\sum_{i=n-q+1}^{n}  (-1)^{i+1} f(2n-i+q+1,i).
\end{equation}
Computing the sums we obtain
\begin{equation}
 \mathrm{res} (C_q)= \frac{n^2+1}{2},
\end{equation}
confirming that the residuum is the same for every column of $M$.

What is now left is checking the diagonals of $M$. The main diagonal of $M$ corresponds to
the middle column of $A$, or $m$-th column,
$\mathrm{diag\,} M = \{f(i,m)\}_{i=0}^n$. It is sorted in increasing order, thus
\begin{equation} \label{resdiag}
 \mathrm{res} \left( \mathrm{diag\,} M \right)
=\sum  _{i=0}^n (-1)^i f(i,m). 
\end{equation}
 The antidiagonal of $M$ corresponds to the $m$-th column of $A$, 
$\mathrm{adiag\,} M = \{f(m,i)\}_{i=0}^n$
and it is 
sorted in decreasing order, thus
\begin{equation} \label{resadiag}
 \mathrm{res} \left( \mathrm{adiag\,} M \right)
=\sum  _{i=0}^n (-1)^{i+1} f(m,i). 
\end{equation}
Both sums in eq. (\ref{resdiag}) and (\ref{resadiag}) can be easily computed. They are both equal to $(n^2+1)/2$, thus
$M$ is indeed a magic square of subtraction.  \Square
\begin{figure}
\begin{equation*}
\left[ \begin {array}{ccccccc} 28&47&16&37&12&35&4
\\ \noalign{\medskip}5&27&48&17&36&11&29\\ \noalign{\medskip}30&6&26&
49&18&42&10\\ \noalign{\medskip}9&31&7&25&43&19&41
\\ \noalign{\medskip}40&8&32&1&24&44&20\\ \noalign{\medskip}21&39&14&
33&2&23&45\\ \noalign{\medskip}46&15&38&13&34&3&22\end {array}
 \right] 
 \left[ \begin {array}{ccccccccc} 37&78&35&66&21&62&15&46&5
\\ \noalign{\medskip}6&38&79&34&65&22&63&14&54\\ \noalign{\medskip}53&
7&39&80&33&64&23&55&13\\ \noalign{\medskip}12&52&8&40&81&32&72&24&56
\\ \noalign{\medskip}57&11&51&9&41&73&31&71&25\\ \noalign{\medskip}26&
58&10&50&1&42&74&30&70\\ \noalign{\medskip}69&27&59&18&49&2&43&75&29
\\ \noalign{\medskip}28&68&19&60&17&48&3&44&76\\ \noalign{\medskip}77&
36&67&20&61&16&47&4&45\end {array} \right] 
\end{equation*}
\caption{Squares of subtraction of order 7 and 9.}\label{exorder79}
\end{figure}

Proposition 1 produces, for $n=5$, the square shown in Figure~\ref{originalexamples}e. How did Kocha\'nski obtain
the remaining three squares of order 5, shown in Figures~\ref{originalexamples}f,g and h? He does not explain
it in the paper, but it is rather easy to guess, following observations made by \cite{Pawlikowska69}.

Keeping the notation used in the proof of Proposition 1, let us recall that the magic square $M$ is obtained
by rotating the initial matrix $A$ and then relocating outside entries to the interior of the square.
Obviously, if $A$ is transposed before rotation and relocation, the resulting square will still be magic.
Transposition of $A$, therefore,  does not affect the magic property. One can also show that any permutation
of columns of $A$ which keeps the middle column in the same place, performed before the rotation
and relocation,  preserves the magic property of the resulting square $M$ (recall that the middle column
of $A$ becomes the diagonal of $M$, so it cannot be moved). The same applies to any permutation of 
rows of $A$ keeping the middle row in the same place.

The remaining three 5th order squares of  Figure~\ref{originalexamples} can be obtained if certain operations of the aforementioned type  are applied to rows and columns of matrix  $A$ of eq. (\ref{defA}) before it is rotated. These are shown below.

\begin{itemize}
 \item For the square shown in  Figure~\ref{originalexamples}f: transpose $A$,
then apply the permutation $(1,2,3,4,5) \to (4,5,3,1,2)$ to the
rows of $A^T$ and  the permutation  $(1,2,3,4,5) \to (4,1,3,5,2)$
to its columns;

 \item  For the square shown in  Figure~\ref{originalexamples}g: 
apply the permutation $(1,2,3,4,5) \to (2,5,3,1,4)$ to the
rows of $A$ and the permutation  $(1,2,3,4,5) \to (1,4,3,2,5)$
to its columns;

 \item  For the square shown in  Figure~\ref{originalexamples}h: 
apply the same permutation  $(1,2,3,4,5) \to (5,2,3,4,1)$ to both the
rows and  columns of $A$.
\end{itemize}
The above operations will produce  three matrices,
\begin{equation*}
A^{\prime}= \left[ \begin {array}{ccccc} 7&10&8&6&9\\ \noalign{\medskip}4&1&3&5&2
\\ \noalign{\medskip}14&11&13&15&12\\ \noalign{\medskip}24&21&23&25&22
\\ \noalign{\medskip}17&20&18&16&19\end {array} \right], 
\,\,
A^{\prime\prime} =\left[ \begin {array}{ccccc} 22&9&12&19&2\\ \noalign{\medskip}25&6&15
&16&5\\ \noalign{\medskip}23&8&13&18&3\\ \noalign{\medskip}21&10&11&20
&1\\ \noalign{\medskip}24&7&14&17&4\end {array} \right], 
\,\,
A^{\prime\prime\prime} = \left[ \begin {array}{ccccc} 5&16&15&6&25\\ \noalign{\medskip}2&19&12
&9&22\\ \noalign{\medskip}3&18&13&8&23\\ \noalign{\medskip}4&17&14&7&
24\\ \noalign{\medskip}1&20&11&10&21\end {array} \right]. 
\end{equation*}
By rotating each one of them and relocating outside entries to the interior we obtain the magic squares of subtraction
shown respectively in Figures \ref{originalexamples}f,
\ref{originalexamples}g, and \ref{originalexamples}h.

Obviously, since Proposition 1 is valid for any odd $n$, we can use it to construct squares of
subtraction of higher orders. Figure~\ref{exorder79} show squares of subtraction of order 7 and 9 obtained by
this method.

\section{Even order}
Just as in the case of magic squares of summation, magic squares of subtraction
of even order are  more difficult to construct. Kocha\'nski gave four examples of squares of order 4, the smallest order
for which  magic squares of subtraction exist,  and these are shown in
Figures~ \ref{originalexamples}a--d. He did not explain how they were constructed. Close inspection of them  reveals 
certain regularities in the arrangement of numbers, and from this  one can guess how the construction probably proceeded.

\begin{figure}
\includegraphics[width=5cm,valign=c]{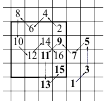}
\,\,\,\,\,\,\,\,$\longrightarrow$\,\,\,\,\,\,\,
{\Large \begin{tabular}{|c|c|c|c|}\hline
1 & 6 & 13 & 2 \\ \hline
10 & 5 &14 & 9  \\ \hline
7 & 12 &11  & 16 \\ \hline
8 & 3 & 4 & 15 \\ \hline
\end{tabular} }
\caption{Construction of the first 4-th order square of subtraction of  Figure~\ref{originalexamples}. Numbers outside
of the framed square are relocated to the interior of the square, and this produces the square on the right.} \label{order4const}
\end{figure}

Most likely,  his first observation was that that the arrangement of odd and even numbers must follow some
pattern reflecting constrains required for the square to be magic. The simplest of such patterns is shown below.
\begin{equation*}
\begin{array}{|c|c|c|c|} \hline 
\mathrm{odd} & \mathrm{even} & \mathrm{odd} & \mathrm{even} \\ \hline 
\mathrm{even} & \mathrm{odd} & \mathrm{even} & \mathrm{odd} \\ \hline
\mathrm{odd} & \mathrm{even} & \mathrm{odd} & \mathrm{even} \\ \hline
\mathrm{even} & \mathrm{odd} & \mathrm{even} & \mathrm{odd} \\ \hline
\end{array} 
\end{equation*}
It satisfies the obvious requirement that the residua of all rows, columns and diagonals
must have the same parity, and that the total number of odd entries must be the same as the total number of 
even entries. He then probably proceeded, by trial and error, to fill this pattern with consecutive odd and even integers, and one of the most natural ways to do this would be to follow a ``zigzag'' path. Such a path does not need to start inside  the 
square, just as in the case of squares of order 5. After some tinkering with such paths, one discovers
that the placement of odd and even numbers shown in Figure~\ref{order4const} fits the bill.
This requires, as in the case of order 5 square, the relocation of numbers which are outside  the square to its interior,
using ``periodic boundary conditions'', that is, left-right and bottom-top wrapping. One thus discovers the square
shown in Figure~ \ref{originalexamples}a. 

\begin{figure}
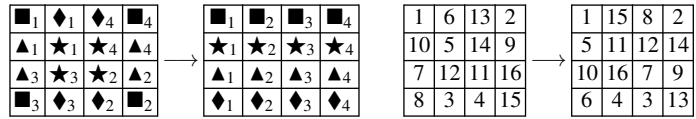

$$
\begin{tabular}{|c|c|c|c|}\hline
$\blacksquare_1$ & $\blacklozenge_1$ & $\blacklozenge_4$ & $\blacksquare_4$ \\ \hline
$\blacktriangle_1$ & $\bigstar_1$ &$\bigstar_4$ & $\blacktriangle_4$  \\ \hline
$\blacktriangle_3$ & $\bigstar_3$ &$\bigstar_2$  & $\blacktriangle_2$ \\ \hline
$\blacksquare_3$ & $\blacklozenge_3$ & $\blacklozenge_2$ & $\blacksquare_2$ \\ \hline
\end{tabular}
\longrightarrow
\begin{tabular}{|c|c|c|c|}\hline
$\blacksquare_1$ & $\blacksquare_2$ & $\blacksquare_3$ & $\blacksquare_4$ \\ \hline
$\bigstar_1$ & $\bigstar_2$ & $\bigstar_3$ & $\bigstar_4$ \\ \hline
$\blacktriangle_1$ & $\blacktriangle_2$ & $\blacktriangle_3$ & $\blacktriangle_4$ \\ \hline
$\blacklozenge_1$ & $\blacklozenge_2$ & $\blacklozenge_3$ & $\blacklozenge_4$ \\ \hline
 \end{tabular}
\,\,\,\,\,\,\,\,\,\,\, \,\,\,\,\,\,\,
\begin{tabular}{|c|c|c|c|}\hline
1 & 6 & 13 & 2 \\ \hline
10 & 5 &14 & 9  \\ \hline
7 & 12 &11  & 16 \\ \hline
8 & 3 & 4 & 15 \\ \hline
\end{tabular}
\longrightarrow
\begin{tabular}{|c|c|c|c|}\hline
1 & 15 & 8 & 2 \\ \hline
5 & 11 &12 & 14  \\ \hline
10 & 16 & 7 & 9 \\ \hline
6 &4  & 3 & 13 \\ \hline
\end{tabular}
$$
\caption{Transformation needed to produce the second magic square of subtraction of Figure~\ref{originalexamples}.}\label{second4}
\end{figure}

The square of   Figure~\ref{originalexamples}b can be produced by
applying to the first square the  transformation shown schematically in Figure~\ref{second4}.
By direct verification one can check that this indeed produces a magic square of subtraction.

Finally, if we apply to both rows and columns of the squares  shown in Figures~ \ref{originalexamples}a and \ref{originalexamples}b the permutation $(1,2,3,4)\to(3,1,4,2)$, we will obtain
the remaining  two squares shown Figures~\ref{originalexamples}c and \ref{originalexamples}d. 
Application of  the same permutation to both  rows and columns 
preserves the magic property, so we are guaranteed that these are  magic squares of subtraction too.

As clever as the above constructions are, they do not seem to be amenable to a generalization to
higher orders. However, Kocha\'nski claims that he found a general method of
construction of doubly-even squares,  of course without revealing any details. In spite of the lack of evidence,
one might  speculate that this general method possibly exploited the existence of 
squares of order 4 to build larger squares  of doubly-even orders. A similar method exists for squares of summation,  and
is sometimes known as the method of composite squares. This is a very
obvious construction, so it seems probable that Kocha\'nski was aware of it. We will now show how this can be done.

Let us assume that $n=4k$, $k \in \NN$. The following proposition constructs magic square of subtraction
of order $n$ if a known square of order 4 is provided.
\begin{proposition}
Let $P$ be a magic square of subtraction of order 4, $J$ be a matrix $4 \times 4$ with all entries equal to 1, 
and let
$Q$ be a block matrix $k \times k$ with entries
$Q_{i,j}= P +\left( k(i-1)+j-1 \right) J $. Then $Q$ is a magic square 
of subtraction of order $4k$.
\end{proposition}
Matrix $Q$ is simply a block matrix with entries consisting of sums of  $P$ and 
 matrices $0,J, 2J\ldots, (k^2-1)J$, added in consecutive rows in increasing order. For example, for $k=3$, it is given by
\begin{equation}
Q=
\left(\begin{array}{lll}
P & P+16J & P+2\cdot 16J\\
P+3\cdot 16J & P+4\cdot 16J & P+5\cdot 16J\\
P+6\cdot 16J & P+7\cdot 16J & P+8\cdot 16J\\
\end{array} \right).    
\end{equation}
Note that each matrix $P+mJ$ is a magic square of subtraction, regardless of $m$,
and its residuum is the same as the residuum of $P$. 

Consider now, for example, the first row of $Q$, which consists of the first
row of $P$ followed by the first row of $P+16J$  and then by the first
row of  $P+32J$, which we will symbolically write as
\begin{equation}
 Q_1= (P_1,  P_1+16 J_1 , P_1+32J_1).
\end{equation}
Note that every entry of $P_1$ is smaller than every entry of $P+16J$,
and every entry of $P+16J$ is smaller than every entry of $P+32J$. For this reason,
\begin{gather}
\mathrm{res}(Q_1)= \mathrm{res}(P_1)+ \mathrm{res}( P_1+16 J_1)+ \mathrm{res}( P_1+32J_1)\\
= \mathrm{res}(P_1)+\mathrm{res}(P_1)+\mathrm{res}(P_1)=3\, \mathrm{res}(P).
\end{gather}
Exactly the same relationship holds for every row and column of $Q$, as well as for
both diagonals, for the same reason as above. Obviously, this observation
also easily  generalizes to arbitrary $k$, when 
the residuum of every row, column and both diagonals of $Q$ would be equal to $k \, \mathrm{res}(P)$,
proving that $Q$ is indeed a magic square of subtraction.$\Square$

For $k=2$, an example of construction of a square of order 8 using this method is shown below. 
\begin{equation*}
P=
 \left[ \begin {array}{cccc} 1&6&13&2\\ \noalign{\medskip}10&5&14&9
\\ \noalign{\medskip}7&12&11&16\\ \noalign{\medskip}8&3&4&15
\end {array} \right],
Q=\left[ \begin {array}{cccc:cccc} 
1&6&13&2&17&22&29&18 \\ 
10&5&14&9&26&21&30&25\\ 
7&12&11&16&23&28&27&32\\ 
8&3&4&15&24&19&20&31\\ \hdashline 
33&38&45&34&49&54&61&50\\  
42&37&46&41&58&53&62&57\\ 
39&44&43&48&55&60&59&64\\ 
40&35&36&47&56&51&52&63
\end {array} \right] 
\end{equation*}
This example uses as $P$ one of the squares of order $4$ constructed by Kocha\'nski, with
residuum 8. The residuum of the resulting square $Q$ is equal to $k \, \mathrm{res}(P)=2\cdot 8 =16$. 

What remains is the case of singly-even squares, of order $n=4k+2$,  $k \in \NN$. Kocha\'nski
apparently tried to construct the smallest order square of subtraction of this type, for
$n=6$, but he failed. For this reason, at the end of his paper, he challenged fellow mathematicians
to find the square of order 6:
\begin{quotation}
\noindent\emph{PROBLEMA I. In Quadrato Senarii, cellulas 36 complectente, numeros progressionis
Arithmeticae, ab 1. ad 36.
 inclusive procedentes, ita disponere, ut subtractionis articio, prius explicato, in omnibus columnis, 
trabibus, \& utraque Diagonali relinquatur numerus 18.}\\
PROBLEM I. In a square of base 6, containing 36 cells, to arrange numbers of arithmetic
sequence, proceeding inclusively from 1 to 36, in such a a way that by the method of subtraction
explained earlier, the number 18 remains in all columns, rows and both diagonals.
\end{quotation}
Although I was not able to discover any general method
of construction of singly-even squares, I found millions of $n=6$ squares by brute-force
computerized search. One  such square is shown below. It has residuum 18, just as Kocha\'nski wished:
\begin{equation*}
\left[
\begin{array}{cccccc} 
1 &11&12&13&29&2 \\ 
17&5 &19&20&8&33 \\ 
34&14&6&9&28&23 \\ \
32&16&10&7&24&25 \\ 
35&26&27&22&21&15 \\ 
3&36&18&31&30&4 
\end{array}\right].
\end{equation*}
The question asked 330 years ago can, therefore, be considered answered, although the answer is not  entirely satisfying. It has been obtained 
with the help of a machine, while Kocha\'nski wanted
to \emph{provocare homines ut vires suas experiantur}, that is, to \emph{provoke people to put to the test their [mental] powers}.
On the other hand, he spent a lot of time and effort designing and attempting to construct computing machines, so perhaps he would
appreciate the fact that such machines were eventually constructed and that they helped to produce the example he asked for.

A general method for construction of singly-even squares of orders higher than 6 remains, for now, an open problem. The known algorithms
for construction of singly even magic squares of summation do not seem to provide any useful insight.

\section{Enumeration}
It is rather surprising to compare the number of existing magic  squares of subtraction
with the number of magic squares of summation. For magic squares of summation, their number 
is known up to order 5 \citep{enum}. It is customary to count not all magic squares, but only ``distinct'' ones,
that is, the number of  equivalence classes with respect to the group of 8 symmetries of the square (dihedral group $D_4$).
Table~\ref{table1}  summarizes currently known enumeration results, counting only the above equivalence classes.
\begin{table}
 \begin{tabular}{l|c|c|c|c|c}
 &$n=2$ & $n=3$ & $n=4$ & $n=5$ & $n=6$  \\ \hline
 number of distinct magic squares  of summation\,\,\,\,\,        &0  & 1 & 880    & 275,305,224 & ?  \\  \hline
 number of distinct magic squares of subtraction & 0 & 0 & 24,488 & ? & ?
\end{tabular}
\caption{Enumeration of magic squares of summation and subtraction.}\label{table1}
\end{table}
In spite of the fact that  magic squares of subtraction are subjected to strongly nonlinear constraints
(equality of residua, which require sorting of vectors), for $n>3$ they are apparently
much more common than magic squares of summation, for which the constraints are linear.  For $n=4$, I
found that the number of distinct magic squares of subtraction is 22,488, which is over 25 times more than
the number of magic squares of summation of the same order. The exact number of squares of subtraction of
order 5 remains unknown, but the non-exhaustive computer search  shows that
their number exceeds $10^9$, possibly greatly, so it is certainly much larger than the number of magic squares of summation of order~5.

Direct computer enumeration of squares of order 5 is probably within the reach of current hardware,  although
I have not been able to complete it yet. Squares of order 6, on the other hand, are probably
outside of the reach of today's computers, and may remain so for a long time (maybe even forever).

\begin{acknowledgement}
 The author acknowledges financial support from the Natural Sciences and
Engineering Research Council of Canada (NSERC) in the form of Discovery Grant.
The computational part of this work was made possible by the facilities of the Shared
Hierarchical Academic Research Computing Network (SHARCNET:
www.sharcnet.ca) and Compute/Calcul Canada.
\end{acknowledgement}

\end{document}